\documentclass[12pt]{article}
\usepackage{amsmath}
\usepackage{amsfonts}

\newtheorem{thm}{{\sc Theorem}}

\newcommand{\qed}{\hspace*{\fill} Q.E.D.}

\title{Addendum To: On fibre space structures of a projective irreducible
       symplectic manifold}
\author{Daisuke Matsushita}
\date{E-mail Address tyler@kurims.kyoto-u.ac.jp}
\begin{document}
\maketitle
%\begin{abstract}
% In this addendum, we prove that every fibre space of a projective
% irreducible manifold is a lagrangian fibration.
%\end{abstract}
%\section{Introduction}
% First we define an {\it irreducible symplectic manifold}.
%
%\noindent
%{\sc Definition. \quad} 
% A simply connected compact K\"{a}hler manifold $X$ is said to
% be {\it irreducible symplectic} if $X$ satisfies
% the following two conditions:
%\begin{enumerate}
% \item There exists non-degenerate holomorphic two form on $X$.
% \item $h^{2,0} = 1$.
%\end{enumerate}
% These manifolds form a building blocks of all Ricci-flat
% K\"{a}hler manifolds together with Calabi-Yau manifolds and
% Abelian varieties due to the Bogomolov decomposition theorem
% \cite{bogomolov}.
% In dimension 2, $K3$ surfaces are the only irreducible symplectic
% manifold, and irreducible symplectic manifolds are considered
% as higher-dimensional analogies of $K3$ surfaces.
% Certain $K3$ surfaces $S$ admit a fibre space structure 
% $f : S \to {\mathbb P}^{1}$ 
% whose general fibre is an elliptic curve. A general fibre of $f$ is
% a lagrangian submanifold
% becauase every fibre is one dimensional.
 In this note, we prove that 
 every fibre space of a projective
 irreducible manifold is a lagrangian fibration.
\begin{thm}
 Let $X$ be a projective irreducible symplectic manifold and
 $f: X \to B$ a fibre space with projective base $B$. 
 Then $f$ is a lagrangian fibration, that is 
 a general fibre of $f$ is 
 a lagrangian submanifold.
\end{thm}
{\sc Remark. \quad} 
 Beauville proves 
 that 
 a Lagrangian fibration is a complete integrable system
 in \cite[Proposition 1]{beauville2}.
%
% $B$3$3$G(B Beauville $B$5$s$NO@J8$r0zMQ$9$k(B
%
% Moreover, he shows that a general fibre of lagrangian fibration is
% a abelian variety.
 Thus, a general fibre of a fibre space of a projective
 irreducible symplectic manifold is an abelian variety.

\noindent
{\sc Remark. \quad}
 Markshevich states in \cite[Remark 3.2]{mark1} that there exists 
 an irreducible 
 symplectic manifold which has a family of 
 non lagrangian tori.
 But this family does not form fibration.
%
% general fibre $B$,(B Abelian variety $B$G$"$k$3$H(B 
%
%
\vspace{5mm}

\noindent
{\sc Proof of Theorem. \quad}
 Let $\omega$ be a nondegenerate two form on $X$ and $\bar{\omega}$
 a conjugate. Assume that $\dim X = 2n$.
 Then $\dim F = n$ \cite[Theorem 2]{matsu}, where $F$ is a
 general fibre of $f$.
 In order to prove that $F$ is a Lagrangian submanifold,
 it is enough to show
$$
 \int_{F} \omega \wedge \bar{\omega} A^{n-2} = 0
$$
 where $A$ is an ample divisor on $X$.
 Let $H'$ be an ample divisor on $B$  and $H := f^{*}H'$.
 Then
$$
 \int_{F} \omega \wedge \bar{\omega} A^{n-2} = 
 c (\omega  \bar{\omega} A^{n-2} H^n ),
$$
 where $c$ is a nonzero constant.
 We shall show $\omega \bar{\omega} A^{n-2} H^n = 0$.
 By \cite[Theorem 4.7]{fujiki}, 
 there exists a bilinear form $q_{X}$
 on $H^{2}(X,{\mathbb C})$ which has the following properties:
$$
 a_0 q_{X}(D,D)^n = D^{2n} \quad D \in H^2(X, {\mathbb C}). 
$$ 
 We consider the following equation,
$$
 a_0 q_{X}(\omega + \bar{\omega} + sA + tH,
   \omega + \bar{\omega} + sA + tH )^n 
 = ( \omega +\bar{\omega} + sA + tH )^{2n}.
$$
 Calculating the left hand side, we obtain
$$
 a_0 (q_{X}(\omega + \bar{\omega}) + s^2 q_{X}(A)
      + 2sq_{X}(\omega + \bar{\omega},A) + 
        2tq_{X}(\omega + \bar{\omega},H)  +
        2stq_{X}(A,H))^n . 
$$
 Since $\omega \in H^{0}(X, \Omega_{X}^{2})$ and
       $A,H \in H^{1}(X,\Omega_{X}^{1})$,
$$
 q_{X}(\omega + \bar{\omega}, A) = 
 q_{X}(\omega + \bar{\omega}, H) = 0
$$
 by \cite[ Th\'{e}or\`{e}me 5]{beauville}.
 Thus
 we can conclude that  $\omega \bar{\omega} A^{n-2} H^n = 0$ by
 comparing $s^{n-2}t^n$ term of both hands sides.
\qed

\begin{flushleft}
 Resarch Institute of Mathematical Science, \\
 Kyoto University. \\ 
 KITASHIRAKAWA, OIWAKE-CHO, KYOTO, 606-01, JAPAN.
\end{flushleft}
\end{document}